\documentclass[12pt]{amsart}
\usepackage{amsthm,amsmath,amssymb,amscd,amsxtra,graphics,psfig}

{\end{list}}
{
   \newtheorem{theorem}{Theorem}[section]                     
        
   \newtheorem{lemma}[theorem]{Lemma}

}
{\theoremstyle{definition}

   \newtheorem{remark}[theorem]{Remark}

}

\newcommand{\QQ}{{\mathbb{Q}}}

\newcommand{\ZZ}{{\mathbb{Z}}}

\newcommand{\cA}{{\mathcal A}}

\newcommand{\cT}{{\mathcal T}}

\newcommand{\bone}{{\boldsymbol{1}}}

\newcommand{\Hom}{{\operatorname{Hom}}}

\newcommand{\eff}{{\operatorname{eff}}}

\newcommand{\dar}{\downarrow}

\newcommand{\isomto}{\stackrel{\sim}{\to}}
\newcommand{\isom}{\simeq}

\newcommand{\setmin}{\,\protect%
\begin{picture}(8,10)\qbezier(1,5.5)(4,4.)(7,2.5)\end{picture}\,}

\begin{document}
\title[Toric residue mirror conjecture]{Toric Residue Mirror
  Conjecture for Calabi-Yau complete intersections}  

\author{Kalle Karu}
\address{Mathematics Department\\ University of British Columbia \\
  1984 Mathematics Road\\
Vancouver, B.C. Canada V6T 1Z2}
\email{karu@math.ubc.ca}

\maketitle

\addtocounter{section}{-1}

\section{Introduction}

The toric residue mirror conjecture of Batyrev and Materov \cite{BM1,BM2}
expresses a toric  residue as a power series whose coefficients are
certain integrals over moduli spaces. This conjecture for Calabi-Yau
hypersurfaces in Gorenstein toric Fano varieties was proved
independently by Szenes and Vergne 
\cite{SV} and Borisov \cite{B2}. We build on the work of these authors
to generalize the residue mirror map to not necessarily reflexive
polytopes. Using this generalization we prove the toric residue mirror
conjecture for Calabi-Yau complete intersections in Gorenstein toric
Fano varieties \cite{BM2}.

We start by introducing notation and explaining the main idea of the
generalization. We work over the field $K=\QQ$.  Let
$\overline{M}\isom \ZZ^d$, let $\Delta\subset 
\overline{M}_K$ be a $d$-dimensional lattice polytope, and let $\cT$
be a coherent
triangulation of $\Delta$, defined by a convex piecewise linear
integral function on $\Delta$. All lattice points in
$\Delta$ are assumed to be vertices of the simplices in $\cT$. We
place $\Delta$ in $M_K = (\overline{M}\times \ZZ)_K$ as
$\Delta\times\{1\}$ and let $C_\Delta\subset M_K$ be the cone over
$\Delta$ with vertex $0$. Then $\cT$ defines a subdivision of
$C_\Delta$ into a fan $\Sigma$. 

The idea of the toric residue mirror conjecture is to relate the
semigroup ring $S_\Delta = K[C_\Delta\cap M]$ to the cohomology of
the fan $\Sigma$. Let $I_\Delta\subset S_\Delta$ be the ideal
generated by monomials $t^m$ where $m\in M$ lies in the interior of
$C_\Delta$. Given general elements $f_0,\ldots, f_d \in S_\Delta^1$
(the superscript denotes the degree), we can construct the toric
residue map \cite{C}:
\[ Res_{(f_0,\ldots,f_d)}:
(I_\Delta/(f_0,\ldots,f_d)I_\Delta)^{d+1} \isomto K. \]
Following \cite{BM1}, we choose a special set of $f_i$ constructed
from a single $f\in S_\Delta$ by partial differentiation.

On the cohomology side, the Poincar\'e dual of the cohomology
$H(\Sigma)$ is the cohomology with compact support
$H(\Sigma,\partial\Sigma)$ \cite{BBFK}. In the top degree we have the
evaluation map 
\[ \langle \cdot \rangle_\Sigma: H^{d+1}(\Sigma,\partial\Sigma)
\isomto K.\]
The residue mirror map takes $I_\Delta^{d+1}$ into
$H^{d+1}(\Sigma,\partial\Sigma)$ so that composition with the
evaluation map gives the toric residue.

The toric residue mirror conjecture of Batyrev and Materov
\cite{BM1,BM2} is a special case of the above formulation. If $\Delta$
is reflexive, it has only one lattice point $0$ in its
interior. Assume that every maximal simplex in $\cT$ has $0$ as a
vertex. Then the projection $q: M_K \to \overline{M}_K$ maps the
fan $\Sigma$ to a complete fan $\overline{\Sigma}$ in
$\overline{M}_K$. (Geometrically, the toric variety of $\Sigma$ is
the total space of a line bundle over the toric variety of
$\overline{\Sigma}$.) The cohomology spaces of the two fans are
isomorphic, hence we can express the toric residue in terms of the
cohomology of $\overline{\Sigma}$.   

In the complete intersection case we use the Cayley
trick \cite{BM2} to construct a polytope $\tilde\Delta \subset
M_K = (\overline{M}\times \ZZ^r)_K$ and a fan $\Sigma$ subdividing
$C_{\tilde\Delta}$. The projection $q: M_K \to \overline{M}_K$ again maps
$\Sigma$ to a complete fan $\overline{\Sigma}$. (The geometry here is
  that the toric variety of $\Sigma$ is the total space of a rank $r$
  vector bundle over the toric variety of $\overline{\Sigma}$.) Thus,
  we can express the toric residue in terms of the cohomology of
  $\overline{\Sigma}$. 

In the complete intersection case the ring $S_{\tilde\Delta}$ is graded by
$\ZZ_{\geq 0}^r$. Restricting the toric residue to a homogeneous component of
$I_{\tilde\Delta}$ defines the mixed  toric residue. We also prove a
conjecture in 
\cite{BM2} relating the mixed residues with mixed volumes of
polytopes.

In the proofs we follow the algebraic approach of Borisov \cite{B2},
but we replace the higher Stanley-Reisner rings with Jeffrey-Kirwan
residues as in \cite{SV}.

{\em Notation.} Given a lattice $M\isom \ZZ^d$, we denote $M_K =
M\otimes K$ and the dual lattice $N=M^* = Hom(M,\ZZ)$. For $u\in M$ and
$w\in N$, we let the pairing be $(w,u)\in K$. Given a homomorphism
$q: M\to M'$ of lattices, we denote the scalar extension $M_K\to
M_K'$ also by $q$.

\section{Cohomology}

We recall the equivariant definition of the cohomology of $\Sigma$
(which is the cohomology of the associated toric variety) \cite{B, BBFK}.

Let $\cA(\Sigma)$ be the ring of $K$-valued conewise polynomial
functions on $\Sigma$, graded by degree. The cohomology
$H(\Sigma)$ is defined as the quotient $\cA(\Sigma)/I$, where $I$ is
the ideal generated by global linear functions. 

One can recover the
Stanley-Reisner description of cohomology as follows. Let 
$v_1,\ldots,v_n$ be the primitive generators of $\Sigma$ (the first
lattice points on the $1$-dimensional cones of $\Sigma$), and let
$\chi_i \in \cA^1(\Sigma)$ be the conewise linear functions defined by 
\[ \chi_i(v_j) = \delta_{i j},\]
where $\delta_{ij}$ is the Kronecker delta symbol. Then $\chi_i$ for
$i=1,\ldots,n$ generate the ring $\cA(\Sigma)$, with relations
generated by monomials $\prod_{i\in I} \chi_i$, where $\{v_i\}_{i\in
  I}$ do not lie in one cone of $\Sigma$. To obtain the cohomology, we
add the linear relations 
\[ \sum_{i=1}^n (w,v_i) \chi_i = 0 \]
for all $w\in N = \Hom(M,\ZZ)$.

Let $\cA(\Sigma,\partial\Sigma)$ be the ideal in $\cA(\Sigma)$ of
functions vanishing on the boundary of $\Sigma$, and let
$H(\Sigma,\partial\Sigma)$ be the quotient
$\cA(\Sigma,\partial\Sigma)/I\cA(\Sigma,\partial\Sigma)$, where $I$ is
the ideal above. It is proved in
\cite{BBFK} that multiplication of functions induces a non-degenerate
bilinear pairing
\[ H^k(\Sigma) \times H^{d+1-k} (\Sigma,\partial\Sigma) \to  H^{d+1}
(\Sigma,\partial\Sigma) \isom K.\]

The isomorphism $H^{d+1}(\Sigma,\partial\Sigma) \isom K$ can be
defined as follows \cite{B}. For $\sigma \in\Sigma$ a maximal cone,
define $\Phi_\sigma = \prod_{v_i\in\sigma} \chi_i|_\sigma$, where
$|_\sigma$ means that we consider $\Phi_\sigma$ as a global polynomial
function on $M_K$ whose restriction to $\sigma$ is the product of
$\chi_i$. Let $Vol(\sigma)$ be the volume of the parallelotype
generated by $v_i \in\sigma$. Equivalently, it is the index of the
lattice generated by $v_i\in \sigma$ in $M$. Now if
$f\in\cA^{d+1}(\Sigma,\partial\Sigma)$, consider the rational function 
\[ \langle f\rangle_\Sigma = \sum_{\sigma \in \Sigma^{d+1}}
\frac{f|_\sigma}{\Phi_\sigma Vol(\sigma)}.\]
By Brion \cite{B} the poles of this rational function cancel out, so
that $\langle f\rangle_\Sigma$ is a constant, thus defining an
isomorphism 
\[ \langle \cdot\rangle_\Sigma: H^{d+1}(\Sigma,\partial\Sigma) \isomto
K.\]

We wish to give another description of the evaluation map using
Jeffrey-Kirwan residues \cite{BV, SV}. The method works best for
complete fans, so let us choose a completion $\hat{\Sigma}$ of $\Sigma$
by adding a ray $K_{\geq 0} v_0$ for some $v_0\in M$ such that $-v_0$ lies in
the interior of $C_\Delta$: 
\[ \hat{\Sigma} = \Sigma \cup \{ K_{\geq 0} v_0+\tau| \tau\in\partial\Sigma\}.\]
We have an embedding $H(\Sigma,\partial\Sigma) \subset H(\hat\Sigma)$
defined by extending a function $f\in \cA(\Sigma,\partial\Sigma)$ by
zero outside the support of $\Sigma$. The evaluation map on $H(\hat\Sigma)$
induces the evaluation map on $H(\Sigma,\partial\Sigma)$.

Let $\hat\pi: \ZZ^{n+1} \to M$ be the $\ZZ$-linear map $e_i\mapsto v_i$ for
$e_0,\ldots, e_n$ the standard basis of $\ZZ^{n+1}$.  The kernel of
$\hat\pi$ is 
$R(\hat\Sigma)$, the group of relations among $v_i$. We also let $x_i$
for $i=0,\ldots,n$ be the standard coordinate functions on
$K^{n+1}$. Given a polynomial function $f(x_0,\ldots,x_n)$, we will
consider its restriction to $R(\hat\Sigma)_K \subset K^{n+1}$. 
   
Let $Q$ be the vector space of $K$-valued rational functions on
$R(\hat\Sigma)_K$ with poles lying along the hyperplanes defined by
$x_i=0$. Any element $g\in Q$ of degree $-(n-d) = -\dim
R(\hat\Sigma)_K$ can be written as a linear combination of basic
fractions $(\prod_{i\in I} x^i)^{-1}$, where the images of $\{x_i\}_{i
  \in I}$ form a basis of the dual vector space
$R(\hat{\Sigma})_K^*$, and degenerate fractions where the linear
forms in the denominator do not span the dual. 

The Jeffrey-Kirwan residue according to Brion and Vergne \cite{BV,SV} is
a linear map
\[ \langle \cdot\rangle_{JK(\hat\Sigma)}:  Q^{-(n-d)} \to K, \]
defined on the degenerate fractions to be zero and on the basic
fractions:
 \[ \langle \frac{1}{\prod_{i\in I} x^i} \rangle_{JK(\hat\Sigma)} =
 \begin{cases} \frac{1}{Vol(\sigma)} &  \text{ if $\{v_i\}_{i\notin I}$
     generate a cone $\sigma\in \hat{\Sigma}$,} \\
0 & \text{otherwise.}
\end{cases}
\]

The evaluation map on $H^{d+1}(\hat\Sigma)$ can be given in terms of
the Jeffrey-Kirwan residue as follows. Let $f(x_0,\ldots,x_n)$ be a
homogeneous polynomial of degree $d+1$. Then
\[ \langle f(\chi_0,\ldots,\chi_n) \rangle_{\hat\Sigma} = \langle
\frac{f(x_0,\ldots,x_n)}{x^\bone} \rangle_{JK(\hat\Sigma)},\]
where $x^\bone = x_0 x_1\cdots x_n$.

\begin{lemma}\label{lem-JK-zero} Let $x^m = x_0^{m_0}\cdots x_n^{m_n} \in K[x_0^{\pm 1},\ldots,
    x_n^{\pm 1}]$ be a monomial of degree $-(n-d)$. If
    $\{v_i\}_{m_i\geq 0}$ do not lie in one cone
    $\sigma\in\hat{\Sigma}$ then 
\[ \langle x^m \rangle_{JK(\hat\Sigma)} = 0.\]
\end{lemma}

{\bf Proof.} Write $x^m = x^{m^+}/x^{m^-}$, where $m^+ = \max(m,0)$
and $m^- = \max(-m,0)$. Then $x^m$ can be expressed as a linear
combination of degenerate fractions and basic fractions of the form
$x^{-l}$, where $0\leq l_i\leq m^-_i$ for $i=0,\ldots,n$. If $\{v_i\}_{m_i\geq
  0}$ do not lie in one cone then for no such $l$ can $\{v_j\}_{l_j =
  0}$ generate a cone in $\hat\Sigma$. \qed

Szenes and Vergne \cite{SV} expressed the previous lemma in terms of the Mori
cone as follows. Call an element $L\in \cA^1(\Sigma)$ ample if it is strictly
convex, and a fan quasi-projective if there exists an ample element. By
the assumption that the triangulation $\cT$ is coherent, the fan
$\Sigma$ is quasi-projective. The classes of ample elements form an open
set in $H^1(\Sigma)$ whose closure is called the ample cone. The dual
of the ample cone in $H^1(\Sigma)^* = R(\Sigma)_K$
is the Mori cone of $\Sigma$. Here $R(\Sigma) = \ker(\pi: \ZZ^n\to M)$,
$\pi(e_i) = v_i$ for $e_1,\ldots,e_n$ the standard basis of $\ZZ^n$.
We denote the lattice points in the Mori cone by $R(\Sigma)_\eff$. 

For the following we need to observe that if $\beta=
(\beta_1,\ldots,\beta_n)\in R(\Sigma)$ is such that
$\{v_i\}_{\beta_i<0}$ lie in 
one cone $\sigma\in \Sigma$, then   $\beta \in R(\Sigma)_\eff$. Indeed,
any ample $L$ can be modified by a global linear function so that it
vanishes on $\sigma$ and is strictly positive outside of $\sigma$,
hence its pairing with $\beta$ is non-negative. 

 \begin{lemma}\label{lem-JK-eff} Let $x^m \in K[x_1,\ldots,
    x_n]$ be a monomial of degree $d+1$, and let $\beta\in R(\Sigma)$.
 If $\beta\notin R(\Sigma)_\eff$ then 
\[ \langle \frac{x^{m-\beta}}{x^\bone} \rangle_{JK(\hat\Sigma)} = 0.\]
\end{lemma}

{\bf Proof.} Since $m_i\geq 0$, we have 
\[   \{i\}_{\beta_i< 0} \subset  \{i\}_{m_i-\beta_i-1\geq 0}.\]
It follows from the previous lemma that the Jeffrey-Kirwan residue is
nonzero only if $\{v_i\}_{\beta_i<0}$ is a subset of a cone
$\sigma\in\Sigma$, hence $\beta \in R(\Sigma)_\eff$. \qed

\section{Toric residues}

We recall the definition of toric residues \cite{C,CDS,BM1}. 

Recall that we defined $S_\Delta$ to be the semigroup ring of
$C_\Delta\cap M$ and $I_\Delta\subset S_\Delta$ the ideal generated by
monomials $t^m$ where $m$ lies in the interior of $C_\Delta$. The
ring $S_\Delta$ is Cohen-Macaulay with dualizing module
$I_\Delta$. Given a regular sequence $f_0,\ldots,f_d \in S_\Delta^1$, the
quotient $S_\Delta/(f_0,\ldots,f_d)$ is again Cohen-Macaulay with
dualizing module $I_\Delta/(f_0,\ldots,f_d) I_\Delta$. It follows that
there exists an isomorphism 
\[   (I_\Delta/(f_0,\ldots,f_d) I_\Delta))^{d+1} \isomto K.\]
This isomorphism, normalized so that the Jacobian of $f_0,\ldots,f_d$
maps to $Vol(\Delta)$ is called the toric residue
$Res_{(f_0,\ldots,f_d)}$. Here $Vol(\Delta)$ is $d!$ times the
$d$-dimensional volume of $\Delta$ ($Vol(\Delta) =
\sum_{\sigma\in\Sigma} Vol(\sigma)$). The Jacobian is defined by
choosing a basis $u_i$ for $M$, letting $t_i = t^{u_i}$, and
considering $S_\Delta \subset K[t_0^{\pm 1},\ldots, t_d^{\pm
    1}]$. Then
\[ Jac_{(f_0,\ldots,f_d)} = \det (t_j \frac{\partial f_i}{\partial
  t_j})_{i,j}.\]
The Jacobian lies in $I_\Delta$ and it does not depend on the chosen
basis.

Following Batyrev and Materov \cite{BM1}, we consider a regular
sequence $f_0,\ldots,f_d$,  where
\[ f_i = t_i \frac{\partial f}{\partial t_i}, \quad i=0,\ldots, d\]
and
\[ f = \sum_{i=1}^n a_i t^{v_i},\]
with $a_i$ parameters in $K$. The Jacobian now becomes the Hessian of
$f$:
\[ H_f = \det( t_i \frac{\partial}{\partial t_i}t_j
\frac{\partial}{\partial t_j} f )_{i,j=0,\ldots,d}.\]
Since $t_i \frac{\partial}{\partial t_i} t^{v_k} = (w_i,v_k) t^{v_k}$,
where $w_0,\ldots,w_d$ is the basis of $N$ dual to
$u_0,\ldots,u_d$, we can write the Hessian as 
\[ H_f  = \det(\sum_{k=1}^n (w_i,v_k) (w_j,v_k) a_k
t^{v_k})_{i,j=0,\ldots,d}.\]
By \cite{CDS} the Hessian can also be expanded as
\[ H_f = \sum_{J\subset\{1,\ldots,n\}; |J| = d+1} V(J)^2 \prod_{i\in
    J} a_i t^{v_i},\]
where $V(J)$ is the volume of the cone generated by $\{v_i\}_{i\in
  J}$ (note that this cone may not be a cone in $\Sigma$). Since
$V(J)\neq 0$ only if $\sum_{i\in J} v_i$ lies in the interior of
$C_\Delta$, it follows that $H_f\in I^{d+1}_\Delta$. When
$f_0,\ldots,f_d$ forms a regular sequence, the Hessian 
$H_f$ does not lie in $(f_0,\ldots,f_d)I_\Delta$, hence the
normalization $Res_{a_1,\ldots, a_n} (H_f) = Vol(\Delta)$ defines a
unique linear map 
\[  Res_{a_1,\ldots, a_n}: (I_\Delta/(f_0,\ldots,f_d)
I_\Delta))^{d+1} \isomto K.\]

\section{The residue mirror map}
 
Let $\pi: \ZZ^n \to M$ be the $\ZZ$-linear map $e_i\mapsto v_i$ for
$i=1,\ldots, n$. We define the residue mirror map on monomials $t^l
\in I^{d+1}_\Delta$ by
\[ RM: t^l \mapsto \sum_{m\in \pi^{-1}(l)} \langle \genfrac(){}{}{x}{a}^m
\frac{1}{x^\bone} \rangle_{JK(\hat{\Sigma})} \]
and extend linearly. Here $x^\bone = x_0 x_1\cdots x_n$, 
\[ \genfrac(){}{}{x}{a}^m = \prod_{i=1}^n \genfrac(){}{}{x_i}{a_i}^{m_i},\]
and the sum on the right hand side is considered as a formal sum over
Laurent monomials in $a_i$. Note that such sums do not form a ring,
however multiplication of a formal sum with a Laurent polynomial in 
$a_i$ is well-defined.

If $l = \pi(m_0)$ for some $m_0\in \ZZ^n_{\geq 0}$, then using
Lemma~\ref{lem-JK-eff}, we have
\[ RM: t^l \mapsto \sum_{\beta\in R(\Sigma)} \langle
\genfrac(){}{}{x}{a}^{m_0 -\beta}
\frac{1}{x^\bone} \rangle_{JK(\hat{\Sigma})} = \sum_{\beta\in
  R(\Sigma)_\eff} \langle \genfrac(){}{}{x}{a}^{m_0 -\beta}
\frac{1}{x^\bone} \rangle_{JK(\hat{\Sigma})}.\]
Here the formal sum is a Laurent series in $a_i$ with support lying in the
cone $-m_0 +R(\Sigma)_\eff$. We denote by $K[[a_1,\ldots,a_n]]$ the
ring of such Laurent series (over all $m_0\in \ZZ^n$).

The following two lemmas and their proofs are only slight
modifications of the ones in \cite{B2}.

\begin{lemma} The map $RM$ takes the subspace $((f_0,\ldots,f_d) I_\Delta)^{d+1}$
  to zero. 
\end{lemma}

{\bf Proof.} Consider the linear map from $S_\Delta$ to the
space of formal sums defined on monomials
\[ t^l \mapsto \sum_{m\in \pi^{-1}(l)} \genfrac(){}{}{x}{a}^m.\]
This is a map of $K[x_1,\ldots,x_n]$ modules if we let $x_i$ act on
$S_\Delta$ by multiplication with $a_i t^{v_i}$, and on the formal sums
by multiplication with $x_i$. 

A linear combination $g$ of $f_0,\ldots, f_d$ is given by
\[ g = \sum_{i=1}^n (w,v_i) a_i t^{v_i} \]
for some $w\in N_K$. Thus, multiplication with $g$ in $S_\Delta$
corresponds to multiplication with $\sum_{i=1}^n (w,v_i) x_i$ in the module
of formal sums. Now $R(\hat{\Sigma})_K \subset K^{n+1}$ is defined by
linear equations
\[ \sum_{i=1}^n (w,v_i) x_i + (w,v_0) x_0 = 0.\]
Hence it suffices to show that 
\[ \langle x_0 \genfrac(){}{}{x}{a}^m \frac{1}{x^\bone}
\rangle_{JK(\hat{\Sigma})} = 0 \]
for any $m\in\ZZ^n$ such that $\pi(m)$ lies in the interior of
$C_\Delta$. By Lemma~\ref{lem-JK-zero}, this residue is nonzero only
if $\{v_0\} \cup \{v_i\}_{m_i\geq 0}$ lie in a single cone of
$\hat\Sigma$; in other words,  $\{v_i\}_{m_i\geq 0}$ lie in a cone on
the boundary of $C_\Delta$. Since $\pi(m)\in Int(C_\Delta)$, this
cannot happen. \qed 

For later use we generalize the situation slightly. Let 
\[ f_\gamma = \sum_{i=1}^n a_i \gamma_i t^{v_i},\]
where $\gamma_i >0$ are defined by a $w_\gamma\in N_K$:
\[ (w_\gamma,v_i) = \frac{1}{\gamma_i}, \quad  i=1,\ldots,n.\]
Let $H_{f_\gamma}$ be the Hessian of $f_\gamma$, and consider the
residue mirror map $RM$ applied to $H_{f_\gamma}$ (the map $RM$ is not
changed by $\gamma$).

\begin{lemma}\label{lem-hess}  We have
\[ RM(H_{f_\gamma}) = \sum_{\sigma\in\Sigma^{d+1}} Vol(\sigma)
\prod_{v_i\in\sigma} \gamma_i.\]
In particular, when $\gamma = \bone$, 
\[ RM(H_f) = Vol(\Delta). \]
\end{lemma}

{\bf Proof.} We follow closely the proof of Borisov \cite{B2}.

The Hessian $H_{f_\gamma}$ has an expression
\[ H_{f_\gamma} = \sum_{J\subset\{1,\ldots,n\}; |J| = d+1} V(J)^2 \prod_{i\in
    J} a_i \gamma_i t^{v_i}.\] 
We lift $v_i$ to $e_i\in \ZZ^n$, then
\[ RM(H_{f_\gamma}) = \sum_{J\subset\{1,\ldots,n\}; |J| = d+1} V(J)^2
  \sum_{\beta\in R(\Sigma)_\eff} \langle x^J \gamma^J
  \genfrac(){}{}{x}{a}^{-\beta} 
\frac{1}{x^\bone} \rangle_{JK(\hat{\Sigma})},\]
where we write $x^J = \prod_{i\in J} x_i$ and similarly for
$\gamma^J$. When $\beta=0$, we have  
 \[ \langle \frac{x^J}{x^\bone} \rangle_{JK(\hat\Sigma)} =
 \begin{cases} \frac{1}{V(J)} & \text{ if $\{v_i\}_{i\in J}$
     generate a cone $\sigma\in \Sigma$}, \\
0 & \text{otherwise.}
\end{cases}
\]
It follows that the contribution from $\beta=0$ to $RM(H_{f_\gamma})$ is
\[ \sum_{\sigma\in\Sigma} Vol(\sigma) \prod_{v_i\in\sigma} \gamma_i, \]
 and it remains to show that the contribution from any $\beta\neq 0$
 is zero. 

Fix $\beta\neq 0$ and consider 
\[ H_{f_\gamma}  = \det A = \det(\sum_{k=1}^n (w_i,v_k) (w_j,v_k) a_k
\gamma_k t^{v_k})_{i,j=0,\ldots,d},\]
where $w_0,\ldots,w_d$ is a basis of $N$. Since we want to prove the
vanishing of the contribution from $\beta$ to $RM(H_{f_\gamma})$, we
are allowed to change $H_{f_\gamma}$ by a nonzero constant, so we may
assume $\{w_j\}$ to be a basis of $N_K$ instead of $N$.  We choose the
basis so that $w_0 = w_\gamma$ and $(w_j,v_0) = 0$ for
$j=1,\ldots,d$.   
Then the first row of the matrix $A$ with index $i=0$ has $j$th entry
\[ \sum_{k=1}^n (w_j,v_k) a_k t^{v_k}. \]
Note that $\sum_{k=1}^n (w_j,v_k) x_k +(w_j,v_0) x_0$ restricts to
zero on $R(\hat\Sigma)_K$. Since for $j=1,\ldots,d$, $(w_j,v_0) = 0$, 
we may set the entries $A_{0,j}$ for $j\neq 0$ to zero. From the entry
$j=0$ we get a factor of $x_0$. 

Let $A_{0,0}$ be the minor of the matrix $A$ obtained by removing the
first row and the first column. Similarly to the case of $A$, we have:
\[ A_{0,0} = \det(\sum_{k=1}^n (w_i,v_k) (w_j,v_k) a_k \gamma_k
t^{v_k})_{i,j=1,\ldots,d} = \sum_{J\subset\{1,\ldots,n\}; |J| = d}
    V(J)^2 \prod_{i\in J} a_i \gamma_i t^{v_i},\]
where now $V(J)$ is the $d$-dimensional volume of the cone generated
    by $\{v_i\}_{i\in J}$. This volume is computed by projecting from
    $v_0$ and using the volume form determined by the basis
    $w_1,\ldots,w_d$. 

By the above discussion, disregarding the nonzero constants, we have
to show that  
\begin{equation} \label{eq1} \sum_{J\subset\{1,\ldots,n\}; |J| = d}
  V(J)^2 \gamma^J \langle x_0 
  \frac{x^J}{x^{\beta+\bone}} \rangle_{JK(\hat\Sigma)} = 0.
\end{equation}
Here $\beta+\bone = (1,\beta_1+1,\ldots,\beta_n+1)$.
By Lemma~\ref{lem-JK-zero}, the Jeffrey-Kirwan residue in the formula
is zero unless $\{v_i\}_{\beta_i\leq 0}$ lie in a cone on the boundary
of $C_\Delta$. Since $\beta$ defines a relation among $v_i$, it
follows that $\{v_i\}_{\beta_i\neq 0}$ lie in a proper face of
$C_\Delta$. Let $C_0$ be the minimal such face. By the same lemma, it
now also follows that for the 
Jeffrey-Kirwan residue to be nonzero, $\{v_i\}_{i\in J}$ must lie in a
face $C_1$ of $C_\Delta$ containing $C_0$. 

If $V(J) \neq 0$ in the sum~(\ref{eq1}) above then $\{v_i\}_{i\in J}$ lie in at
most one codimension $1$ face $C_1$ of $C$. Let us fix $C_1$ and prove
\begin{equation} \label{eq2}\sum V(J)^2 \gamma^J \langle x_0
  \frac{x^J}{x^{\beta+\bone}} \rangle_{JK(\hat\Sigma)} = 0,
\end{equation}
where the sum now runs over all $J\subset\{1,\ldots,n\}$,  $|J| = d$
such that $\{v_i\}_{i\in J}$ lie in the face $C_1$. We get the sum
(\ref{eq2}) from (\ref{eq1}) by formally
setting $\gamma_i = 0$ for $v_i\notin C_1$, hence going back to the
determinantal form, we can write the sum (\ref{eq2}) as
\[ \langle \det(B) \frac{x_0}{x^{\beta+\bone}}
  \rangle_{JK(\hat\Sigma)}, \]
where $B$ is the matrix
\[ B = (\sum_{v_k\in C_1} (w_i,v_k) (w_j,v_k)
\gamma_k x_k)_{i,j=1,\ldots,d}.\]
Choose $w_1$ so that 
\[ (w_1,v_k) = \frac{1}{\gamma_k}, \quad v_k\in C_1.\]
Such $w_1$ can be taken as a linear combination of $w_\gamma$ and
$w_1' \in N_K$ vanishing on $C_1$. Then the $j$th entry in the first
row of $B$ is
\[  \sum_{v_k\in C_1} (w_j,v_k) x_k.\]
Since $\sum_{k=1}^n (w_j,v_k) x_k$ restricts to zero on
$R(\hat\Sigma)_K$, we may replace the $j$th entry by
\[  - \sum_{v_k\notin C_1} (w_j,v_k) x_k.\]
After doing this replacement, Borisov \cite{B2} showed that the
support of $\det B$ does not intersect any codimension $1$ face of
$C_\Delta$ containing $C_0$, hence the Jeffrey-Kirwan residue above is
zero. Let us recall his argument.

Choose $w_2,\ldots,w_{r+1}$, where $r= d+1 -\dim C_0$, so that they
vanish on $C_0$. (This choice is made independent of the choice of
$C_1$.) Suppose a monomial $x^I$ that occurs in $\det B$ with 
nonzero coefficient is supported in a codimension $1$ face $C_1'$ of
$C_\Delta$ containing $C_0$. Then we can write $I =
\{i_1,\ldots,i_d\}$, where $v_{i_1} \in C_1'\setmin C_1$ and $
v_{i_2},\ldots, v_{i_d} \in C_1'\cap C_1$. Here $x_{i_1}$ comes from
the first row of the matrix $B$ and $x_{i_2},\ldots, x_{i_d}$ from the
rows $2,\ldots,d$. Because $C_1' \neq C_1$, a nontrivial linear
combination of $w_2,\ldots,w_{r+1}$ vanishes on $C_1'\cap C_1$. It follows that
$x_{i_2},\ldots, x_{i_d}$ do not occur in the nonzero minors of $B$ 
constructed from rows $2,\ldots, r+1$.  \qed

Let $P(x_1,\ldots, x_n)\in  K[x_1,\ldots,x_n]$ be a homogeneous
polynomial of 
degree $d+1$ such
that $P(a_1 t^{v_1},\ldots,a_n t^{v_n}) \in I_\Delta$. It is known
that the residue
\[ Res_{a_1,\ldots,a_n} P(a_1 t^{v_1},\ldots,a_n t^{v_n}) \]
is a rational function in $a_i$ with denominator the principal
determinant $E_f$ \cite{CDS}. The support of $E_f$ is the secondary polytope of
$\Delta$, with vertices corresponding to coherent triangulations of
$\Delta$. Consider the vertex corresponding to the triangulation
$\cT$ and expand $Res_{a_1,\ldots,a_n} P(a_1 t^{v_1},\ldots,a_n
t^{v_n})$  in a Laurent series at that vertex. Since
the inner cone to the secondary polytope at the vertex corresponding
to $\cT$ is the cone $R(\Sigma)_\eff$, the expansion of the residue
lies in the ring that we denoted $K[[a_1,\ldots, a_n]]$. We claim that this
expansion is precisely the one given by the residue mirror map
$RM$. Indeed, modulo the ideal $(f_0,\ldots, f_d)$, we can express
\[  P(a_1 t^{v_1},\ldots,a_n t^{v_n}) =
\frac{g(a_1,\ldots,a_n)}{E_f(a_1,\ldots,a_n)} H_{f},\] 
for some polynomial $g(a_1,\ldots,a_n)$. Then
\[ E_f Res_{a_1,\ldots,a_n} P(a_1 t^{v_1},\ldots,a_n t^{v_n}) = g
Res_{a_1,\ldots,a_n} H_f = g Vol(\Delta),\]
and the same formula holds if we replace $Res_{a_1,\ldots,a_n}$ by
$RM$. Since $E_f$ has a unique inverse in  $K[[a_1,\ldots, a_n]]$, we
get that the two Laurent series are equal. We state this as a theorem.

\begin{theorem}\label{thm-1} Let $P(x_1,\ldots, x_n)\in
  K[x_1,\ldots,x_n]$ be a 
  homogeneous polynomial of degree $d+1$ such 
that $P(a_1 t^{v_1},\ldots,a_n t^{v_n}) \in I_\Delta$. The Laurent
series expansion of 
\[ Res_{a_1,\ldots,a_n} P(a_1 t^{v_1},\ldots,a_n t^{v_n})\]
at the vertex of the secondary polytope of $\Delta$
  corresponding to the triangulation $\cT$ is
\[ Res_{a_1,\ldots,a_n} P(a_1 t^{v_1},\ldots,a_n t^{v_n}) =
\sum_{\beta\in R(\Sigma)_\eff} \langle P(x_1,\ldots,x_n) 
\frac{1}{x^{\beta+\bone}}\rangle_{JK(\hat{\Sigma})} a^\beta.\]
\end{theorem} \qed

In particular, the coefficient of $a^\beta$ in the series above does
not depend on the chosen completion $\hat{\Sigma}$ of the fan $\Sigma$.

\section{Morrison-Plesser fans}

Consider one coefficient of the series in Theorem~\ref{thm-1}:
\[ \langle P(x_1,\ldots,x_n)
\frac{1}{x^{\beta+\bone}}\rangle_{JK(\hat{\Sigma})}.\]
Our goal in this section is to construct a new complete fan
$\hat{\Sigma}_\beta$, the Morrison-Plesser fan, such that the
Jeffrey-Kirwan residue above can be identified with the evaluation map
applied to a top degree cohomology class in
$H(\hat{\Sigma}_\beta)$. In the next
sections we apply this construction to other complete projective fans.

It turns out that $\hat{\Sigma}_\beta$ has a natural description in
terms of Gale dual configurations \cite{SV}. We translate these dual
notions into the more conventional setting of fans.

Let us start by recalling the construction of the fan $\hat\Sigma$ as a
quotient, corresponding to the construction of a toric variety as a
GIT quotient. First note that $\hat\Sigma$ is projective. One can
extend a strictly convex conewise linear function on $\Sigma$ to such
a function $L$ on $\hat\Sigma$ by choosing $L(v_0) \gg 0$.  
Consider the exact sequence
\[ 0\to R(\hat \Sigma) \to \ZZ^{n+1} \stackrel{\hat\pi}{\to} M, \]
and fix an ample class $[L]\in H^1(\hat\Sigma) \isom
R(\hat\Sigma)_K^*$. Then the pair $(\hat\pi,[L])$ determines the fan
$\hat\Sigma$ completely as follows. The Gale dual of a cone
$\sigma\subset M_K$ generated by 
$\{v_i\}_{i\in I}$ is the cone in  $R(\hat\Sigma)_K^*$ generated by
the images of $\{e_i^*\}_{i\notin I}$ under the map $(\ZZ^{n+1})^*
\to R(\hat\Sigma)^*$. Then $\sigma\in\hat\Sigma$ if and only if its
Gale dual contains $[L]$ in its interior. The completeness of the fan
$\hat\Sigma$ corresponds to the condition that the images of
$e_0^*,\ldots,e_n^*$ in $R(\hat\Sigma)_K^*$  lie in an open
half-space; $\hat\Sigma$ being simplicial is equivalent to the
condition that $[L]$ does not lie in a smaller dimensional cone
generated by the images of a subset of $e_i^*$.  The Gale dual cones
also determine the Jeffrey-Kirwan 
residue, 
and hence the evaluation map in the cohomology of $\hat\Sigma$. The 
volume of a cone $\sigma\in\hat\Sigma$ is equal to the volume of its
Gale dual if $\hat\pi$ is surjective; otherwise the volumes differ by a
constant factor, the index $[M:\hat\pi(\ZZ^{n+1})]$.

Let us fix $\beta \in \ZZ^{n+1}$ (take $\beta_0 = 0$ if $\beta\in
 R(\Sigma)\subset \ZZ^n$) and write $\beta = \beta^+ - 
 \beta^-$, where $\beta^\pm_i = \max(\pm \beta_i,0)$. Denote 
 $|\beta^+| = \sum_{i=0}^n \beta_i^+$. The following construction of
 $\hat\Sigma_\beta$ only depends on $\beta^+$.

Let  $\rho: \ZZ^{n+1} \to \ZZ^{n+1+|\beta^+|}$ be the product of diagonal
embeddings $\ZZ\to \ZZ^{1+\beta_i^+}$ for $i=0,\ldots,n$. Define
$M_\beta$ as the pushout of $\rho$ and $\hat\pi$:
$$\begin{array}{rcl}
\ZZ^{n+1+|\beta^+|} & \to & M_\beta \\
\uparrow \rho & & \uparrow \\
\ZZ^{n+1} & \stackrel{\hat\pi}{\to} & M.
\end{array}$$
 In other words, 
\[ M_\beta =  (\ZZ^{n+1+|\beta^+|} \times M)/\ZZ^{n+1},\]
where $\ZZ^{n+1}$ is mapped to the product diagonally. Since $\rho$ embeds
$\ZZ^{n+1}$ in $\ZZ^{n+1+|\beta^+|}$ as a direct summand, $M_\beta$ has no
torsion. From the pushout diagram we also get an exact sequence
\[ 0\to R(\hat\Sigma) \to \ZZ^{n+1+|\beta^+|}
\stackrel{\hat\pi_\beta}{\to} M_\beta \]
and an isomorphism between the cokernels of $\hat\pi$ and
$\hat\pi_\beta$.  
Let $\hat\Sigma_\beta$ be the fan defined by the pair
$(\hat\pi_\beta,[L])$. 

We denote the basis of $\ZZ^{n+1+|\beta^+|}$ by
$\{e_{i,j}\}_{i=0,\ldots,n; j=0,\ldots,\beta^+_i}$ and the
corresponding generators of the fan $\hat\Sigma_\beta$ by
$v_{i,j}\in M_\beta$. The images of the dual basis elements 
$e_{i,j}^*$ under $(\ZZ^{n+1+|\beta^+|})^* \to R(\hat\Sigma)^*$
coincide with the images of $e_i^*$ under $(\ZZ^{n+1})^* \to
R(\hat\Sigma)^*$. It follows from this that $\hat\Sigma_\beta$ is
complete and simplicial. Moreover, the Jeffrey-Kirwan residue
in the fan $\hat\Sigma_\beta$ of a rational function in the variables
$x_{i,j}$ is equal to the Jeffrey-Kirwan residue in the fan
$\hat\Sigma$ of the same function but with $x_{i,j}$ replaced by
$x_i$:
\[ \langle f(x_{i,j})\rangle_{JK(\hat{\Sigma}_\beta)} = \langle
f(x_{i})\rangle_{JK(\hat{\Sigma})}.\]

Now consider the Jeffrey-Kirwan residue at the beginning of this
section. We can express it as:
\begin{align*} \langle P(x_1,\ldots,x_n)
\frac{1}{x^{\beta+\bone}}\rangle_{JK(\hat{\Sigma})} & = 
\langle \frac{P(x_{1,0},\ldots,x_{n,0}) x_{1,0}^{\beta^-_1}\cdots
x_{n,0}^{\beta^-_n}}{\prod_{i,j} x_{i,j}}\rangle_{JK(\hat{\Sigma}_\beta)}\\
&= \langle P(\chi_{1,0},\ldots,\chi_{n,0}) \chi_{1,0}^{\beta^-_1}\cdots
\chi_{n,0}^{\beta^-_n} \rangle_{\hat\Sigma_\beta},
\end{align*}
where we have denoted by $\chi_{i,j}$ the generators of the cohomology
of $\hat\Sigma_\beta$ corresponding to $v_{i,j}$. Let us call
\[ \Phi_\beta = [\chi_{1,0}^{\beta^-_1}\cdots \chi_{n,0}^{\beta^-_n}] \in
H(\hat\Sigma_\beta) \]
the Morrison-Plesser class. Then we have:

\begin{theorem}\label{thm-2} Let $P(x_1,\ldots, x_n)\in
  K[x_1,\ldots,x_n]$ be a 
  homogeneous polynomial of degree $d+1$ such 
that $P(a_1 t^{v_1},\ldots,a_n t^{v_n}) \in I_\Delta$. The Laurent
series expansion of
\[ Res_{a_1,\ldots,a_n} P(a_1 t^{v_1},\ldots,a_n t^{v_n})\]
 at the vertex of the secondary polytope of $\Delta$
  corresponding to the 
triangulation $\cT$ is
\[ Res_{a_1,\ldots,a_n} P(a_1 t^{v_1},\ldots,a_n t^{v_n}) =
\sum_{\beta\in R(\Sigma)_\eff} \langle
P(\chi_{1,0},\ldots,\chi_{n,0}) \Phi_\beta \rangle_{\hat\Sigma_\beta}
a^\beta.\] 
\end{theorem} \qed

\begin{remark} Let $\Sigma_\beta$ be the fan obtained from
  $\hat\Sigma_\beta$ by removing the ray generated by $v_{0,0}$ and
  all cones containing it. Similarly to $\Sigma$, the fan
  $\Sigma_\beta$ is a subdivision of a pointed cone in
  $M_{\beta,K}$. The fan $\Sigma_\beta$ does not depend on the
  completion $\hat\Sigma$ and it can be constructed directly from
  $\Sigma$ by a construction similar to $\hat\Sigma_\beta$. It is also
  possible to show (considering $\Phi_\beta \in \cA(\Sigma_\beta)$):
\[ P(\chi_{1,0},\ldots,\chi_{n,0})\Phi_\beta \in \cA(\Sigma_\beta,\partial
  \Sigma_\beta),\]
hence we can write the series in Theorem~\ref{thm-2} as
\[ Res_{a_1,\ldots,a_n} P(a_1 t^{v_1},\ldots,a_n t^{v_n}) =
\sum_{\beta\in R(\Sigma)_\eff} \langle
P(\chi_{1,0},\ldots,\chi_{n,0}) \Phi_\beta \rangle_{\Sigma_\beta}
a^\beta.\] 
This gives an expansion of the residue independent from the completion
$\hat\Sigma$. However, neither $P(\chi_{1,0},\ldots,\chi_{n,0})$ nor
$\Phi_\beta$ 
may vanish on $\partial\Sigma_\beta$, hence we can not consider
$\Phi_\beta$ as an element in $H(\Sigma_\beta)$ or
$H(\Sigma_\beta,\partial\Sigma_\beta)$.
\end{remark}

\section{Calabi-Yau hypersurfaces}

In this section we explain how the toric residue mirror conjecture for
Calabi-Yau hypersurfaces in Gorenstein toric Fano varieties \cite{BM1,
  B2, SV} follows from Theorem~\ref{thm-1}.

Assume that the polytope $\Delta$ is reflexive; that means, its polar
is also a lattice polytope. Then $0\in \Delta$ is the unique lattice
point in the interior of $\Delta$. We assume that $0$ is a vertex of
every maximal simplex in $\cT$. Let the generators of the fan $\Sigma$
in $M = \overline{M}\times \ZZ$ be $v_i= (\bar{v_i},1)$ for
  $i=1,\ldots,n$ and $v_{n+1} = (0,1)$. Also let $q: M \to
  \overline{M}$ be the projection. Then $q$ maps the fan $\Sigma$ to a
  complete fan $\overline\Sigma$ and we have isomorphisms:
\[  H^i(\overline\Sigma) \stackrel{q^*}{\to} H^i(\Sigma)
\stackrel{\chi_{n+1}}{\longrightarrow} H^{i+1} (\Sigma,\partial\Sigma).\]
These isomorphisms are compatible with evaluation: if $P(x_1,\ldots,
x_n)$ is a homogeneous polynomial of degree $d$ then
\[ \langle P(\bar\chi_1,\ldots,\bar\chi_n)\rangle_{\overline\Sigma} =
\langle \chi_{n+1} P(\chi_1,\ldots,\chi_n)\rangle_{\Sigma},\]
where $\bar\chi_i$ are the generators of the cohomology of
$\overline\Sigma$ corresponding to $\bar{v}_i$. We wish to give a
similar correspondence between the Jeffrey-Kirwan residues in
$\overline\Sigma$ and $\hat\Sigma$. 

Let us choose the completion $\hat\Sigma$ by taking $v_0 = (0,-1)$,
and consider the commutative diagram
\[ \begin{array}{ccccccc}
0 & \to & R(\hat\Sigma) & \to & \ZZ^{n+2} & \to & M \\
  &     & \dar          &     & \dar p    &     & \dar q \\
0 & \to & R(\overline\Sigma) & \to & \ZZ^{n} & \to & \overline{M},
\end{array}
\]
where the middle vertical map is defined by $p(e_i) = e_i$ for
$i=1,\ldots,n$ and $p(e_0) = p(e_{n+1}) = 0$. It follows that
functions defined on $R(\hat\Sigma)_K$ by $x_i$ for $i=1,\ldots,n$ are
the pullbacks of functions defined by $x_i$ on
$R(\overline\Sigma)_K$. The hyperplanes defined by $x_0=0$ and $x_{n+1}
= 0$ map onto $R(\overline\Sigma)_K$. Comparing the volumes of cones in
$\overline\Sigma$ and in $\hat\Sigma$, we get
\[ \langle x^m \rangle_{JK(\overline{\Sigma})} = 
 \langle x^m
\frac{1}{x_0}\rangle_{JK(\hat{\Sigma})}\]
for any Laurent monomial $x^m\in K[x_1^{\pm 1},\ldots, x_n^{\pm
    1}]$. If $l>0$ then 
\[  \langle x_0^l x^m \frac{1}{x_0}\rangle_{JK(\hat{\Sigma})} = 0\]
because all linear forms in the denominator are pulled back from
$R(\overline\Sigma)_K^*$, hence they do not span
$R(\hat\Sigma)_K^*$. Using the linear relation $-x_0 + x_1+ \ldots +
x_{n+1} = 0$ on   $R(\hat\Sigma)_K$, we get for $k\geq 0$
\begin{align*} \langle x^m x_{n+1}^k
  \frac{1}{x_0}\rangle_{JK(\hat{\Sigma})} 
&=\langle
  x^m(x_0-x_1-\ldots-x_n)^k\frac{1}{x_0} \rangle_{JK(\hat{\Sigma})} \\
&=\langle
  x^m(-x_1-\ldots-x_n)^k\frac{1}{x_0} \rangle_{JK(\hat{\Sigma})} \\
&=\langle
  x^m(-x_1-\ldots-x_n)^k \rangle_{JK(\overline{\Sigma})}.
\end{align*}

Let $R(\overline\Sigma)$ be the group of relations among
$\bar{v}_i$. We have an isomorphism 
\begin{align*}
 R(\Sigma) & \to R(\overline\Sigma) \\
(\beta_1,\ldots,\beta_{n+1}) & \mapsto (\beta_1,\ldots,\beta_{n}),
\end{align*}
with inverse defined by $\beta_{n+1} = -\beta_1-\ldots - \beta_n$. The
dual map $H^1(\overline\Sigma) \to H^1(\Sigma)$ identifies the ample
cones of the two fans, hence the map above identifies the Mori
cones. Note also that if 
$\beta\in R(\Sigma)_\eff$ then $\beta_{n+1} \leq 0$ because
$-\chi_{n+1}$ is convex and so it lies in the ample cone of $\Sigma$. 

For $\beta\in R(\overline\Sigma)_\eff$, let $\overline\Sigma_\beta$
be the Morrison-Plesser fan constructed from $\overline\Sigma$. Define
the Morrison-Plesser class $\Phi_\beta \in H(\overline\Sigma_\beta)$:
\[ \Phi_\beta = \bar\chi^{\beta^-} = \bar\chi_{1,0}^{\beta^-_1} \cdots
\bar\chi_{n,0}^{\beta^-_n} (-\bar\chi_{1,0} - \ldots
-\bar\chi_{n,0})^{\beta_1+\ldots + \beta_n},\]
where $\bar\chi_{i,j}$ are the generators of the cohomology of
$\overline\Sigma_\beta$ corresponding to $\bar{v}_{i,j}$. 

The ideal $I_\Delta\subset S_\Delta$ is principal, generated by
$t^{v_{n+1}}$. Consider one coefficient in the series of Theorem~\ref{thm-1}
applied to the polynomial $ a_{n+1} t^{v_{n+1}} P(a_1
t^{v_1},\ldots,a_n t^{v_n}) \in I_\Delta$:
\begin{align*}  \langle x_{n+1} P(x_1,\ldots,x_n) \frac{1}{x^{\beta+\bone}}
\rangle_{JK(\hat{\Sigma})}  &=
   \langle  P(x_1,\ldots,x_n) \frac{x^{\beta^-}}{x^{\beta^+}}
   \frac{1}{x_0 x_1\cdots x_n} \rangle_{JK(\hat{\Sigma})} \\ &=
  \langle  P(x_1,\ldots,x_n) \frac{x^{\beta^-}}{x^{\beta^+}}
   \frac{1}{x_1\cdots x_n} \rangle_{JK(\overline{\Sigma})} \\ &=
 \langle P(\bar\chi_{1,0},\ldots,\bar\chi_{n,0}) \Phi_\beta
\rangle_{\overline\Sigma_\beta}.
\end{align*}

Thus, we get:

\begin{theorem}\label{thm-3} Let $P(x_1,\ldots, x_n)\in
  K[x_1,\ldots,x_n]$ be a homogeneous polynomial of degree $d$. The Laurent
series expansion of $Res_{a_1,\ldots,a_n, a_{n+1} = 1} (t^{v_{n+1}}
P(a_1 t^{v_1},\ldots,a_n t^{v_n}))$ at the vertex of the secondary
polytope of $\Delta$ corresponding to the 
triangulation $\cT$ is
\[ Res_{a_1,\ldots,a_n}(t^{v_{n+1}} P(a_1 t^{v_1},\ldots,a_n t^{v_n})) =
\sum_{\beta\in R(\overline\Sigma)_\eff} \langle
P(\bar\chi_{1,0},\ldots,\bar\chi_{n,0}) \Phi_\beta
\rangle_{\overline\Sigma_\beta} a^\beta.\] \qed
\end{theorem}

In \cite{BM1} the parameters $a_i$ differ by a sign from the ones used
here. This introduces a sign difference in the definition of
$\Phi_\beta$ and in the Laurent series expansion.

\section{Complete intersections}

In this section we prove the toric residue mirror conjecture for
Calabi-Yau complete intersections in Gorenstein toric Fano varieties
\cite{BM2}. The construction relies on the Cayley trick \cite{BM2} and
the proof is completely analogous to the hypersurface case.

Let $\Delta\in \overline{M}_K$ be a reflexive polytope ($\Delta =
\nabla^*$ in \cite{BM2}), and $\cT$ a coherent triangulation of
$\Delta$ such that $0\in\Delta$ is a vertex of every maximal
simplex. Let $\overline\Sigma$ be the complete simplicial fan in
$\overline{M}_K$ defined by $\cT$. Denote by
$\{\bar{v}_1,\ldots,\bar{v}_n\}$ the primitive generators of
$\overline\Sigma$, and let $L$ be the conewise linear function on
$\overline\Sigma$ such
that $L(\bar{v}_i) = 1$ for $i=1,\ldots,n$. A {\em nef-partition}
\cite{B1} of $L$ is an expression
\[ L= l_1 +l_2 +\ldots + l_r, \]
where $l_i$ are integral non-negative convex conewise linear functions on
$\overline\Sigma$. We assume that all $l_i\neq 0$. A nef-partition
defines a partition of 
$\{1,\ldots,n\}$ into a disjoint union  $E_1\cup\ldots\cup E_r$, where
$E_j = \{i | l_j(v_i) = 1\}$. Let 
\[ \Delta_j = conv(\{0\}\cup \{ v_i\}_{i\in E_j}).\]

Let $M = \overline{M}\times \ZZ^r$. Define the Cayley polytope 
\[ \tilde\Delta = \Delta_1 *    \cdots *     \Delta_r =
conv(\Delta_1\times \{(0,e_1)\} \cup\ldots \cup \Delta_r\times \{(0,e_r)\}),\]
where $e_1,\ldots,e_r$ is the standard basis of $\ZZ^r$, and let
$C_{\tilde{\Delta}}$ be the cone over $\tilde\Delta$. The lattice
points in $\tilde\Delta$ are $v_i= (\bar{v}_i,e_j)$ for
$i=1,\ldots,n$, where $i\in E_j$ and $v_{n+j} = (0,e_j)$ for
$j=1,\ldots,r$. The triangulation 
$\cT$ defines a triangulation $\tilde\cT$ of $\tilde\Delta$, hence a
simplicial subdivision of the cone $C_{\tilde\Delta}$ into a fan
$\Sigma$ as follows. Let the maximal cones of $\Sigma$ be generated by 
\[ \{v_{n+1},\ldots, v_{n+r}\} \cup \{v_i\}_{\bar{v}_i \in \sigma} \]
for some maximal cone $\sigma\in \overline\Sigma$. 

Let $q: M\to \overline{M}$ be the projection, mapping the fan $\Sigma$
to the fan $\overline\Sigma$. Since every maximal cone in $\Sigma$ is
the product of a cone in $\overline\Sigma$ with the simplicial cone
generated by $\{v_{n+1},\ldots, v_{n+r}\}$, we get isomorphisms
\[  H^i(\overline\Sigma) \stackrel{q^*}{\to} H^i(\Sigma)
\stackrel{\chi_{n+1}\cdots \chi_{n+r}}{\longrightarrow} H^{i+r} (\Sigma,\partial\Sigma).\]
These isomorphisms are compatible with evaluation maps: if $P(x_1,\ldots,
x_n)$ is a homogeneous polynomial of degree $d$ then
\[ \langle P(\bar\chi_1,\ldots,\bar\chi_n)\rangle_{\overline\Sigma} =
\langle \chi_{n+1} \cdots \chi_{n+r}
P(\chi_1,\ldots,\chi_n)\rangle_{\Sigma}.\] 

We complete $\Sigma$ to $\hat\Sigma$ by adding the ray generated by
$v_0 = -v_{n+1}-\ldots -v_{n+r}$ and consider the commutative diagram

\[ \begin{array}{ccccccc}
0 & \to & R(\hat\Sigma) & \to & \ZZ^{n+r+1} & \to & M \\
  &     & \dar          &     & \dar p    &     & \dar q \\
0 & \to & R(\overline\Sigma) & \to & \ZZ^{n} & \to & \overline{M},
\end{array}
\]
where the middle vertical map is defined by $p(e_i) = e_i$ for
$i=1,\ldots,n$ and $p(e_i) = 0$ for $i=0,n+1,\ldots,n+r$. The
functions defined on $R(\hat\Sigma)_K$ by $x_i$ for $i=1,\ldots,n$ are
pullbacks of functions on
$R(\overline\Sigma)_K$; the hyperplanes defined by $x_i=0$ for
$i=0,n+1,\ldots,n+r$  map onto $R(\overline\Sigma)_K$.
One easily checks (for example, using the comparison of the evaluation
maps in $H(\overline\Sigma)$ and $H(\hat\Sigma)$) that
\[ \langle x^m \rangle_{JK(\overline{\Sigma})} = 
 \langle x^m \frac{1}{x_0}\rangle_{JK(\hat{\Sigma})}\]
for any Laurent monomial $x^m\in K[x_1^{\pm 1},\ldots, x_n^{\pm
    1}]$. As in the previous section, we also have
\[  \langle x_0^l x^m \frac{1}{x_0}\rangle_{JK(\hat{\Sigma})} = 0\]
for any $l>0$, and using the relations $-x_0 + x_{n+j} +\sum_{i\in
  E_j} x_i = 0$ on $R(\hat\Sigma)_K$, we get for $k_i\geq0$
\[ \langle x^m x_{n+1}^{k_1} \cdots
  x_{n+r}^{k_r}\rangle_{JK(\hat{\Sigma})} =  
\langle  x^m (-\sum_{i\in E_1} x_i)^{k_1} \cdots (-\sum_{i\in E_r}
x_i)^{k_r}\rangle_{JK(\overline{\Sigma})}. 
\]

Forgetting the last $r$ coordinates of vectors in $\ZZ^{n+r}$, we get
an isomorphism  
\begin{align*}
 R(\Sigma) & \to R(\overline\Sigma) \\
(\beta_1,\ldots,\beta_{n+r}) & \mapsto (\beta_1,\ldots,\beta_{n}),
\end{align*}
with inverse defined by $\beta_{n+j} = -\sum_{i\in E_j} \beta_i$. This
isomorphism identifies the Mori cones of $\Sigma$ and
$\overline\Sigma$. If $\beta\in R(\Sigma)_\eff$ then $\beta_{n+j}
\leq 0$ because $-\chi_{n+j}$ lies in the ample cone of $\Sigma$.

For $\beta\in R(\overline\Sigma)_\eff$, let $\overline\Sigma_\beta$
be the Morrison-Plesser fan constructed from $\overline\Sigma$. Define
the Morrison-Plesser class $\Phi_\beta \in H(\overline\Sigma_\beta)$:
\[ \Phi_\beta = \bar\chi^{\beta^-} = \bar\chi_{1,0}^{\beta^-_1} \cdots
\bar\chi_{n,0}^{\beta^-_n} (-\sum_{i\in
  E_1}\bar\chi_{i,0})^{\sum_{i\in E_1} \beta_i} \cdots (-\sum_{i\in
  E_r}\bar\chi_{i,0})^{\sum_{i\in E_r} \beta_i}.
\]

The ideal $I_{\tilde\Delta} \subset S_{\tilde\Delta}$ is again
principal, generated by 
$t^{v_{n+1}} \cdots t^{v_{n+r}}$. Thus, we have

\begin{theorem}\label{thm-4} Let $P(x_1,\ldots, x_n)\in
  K[x_1,\ldots,x_n]$ be a homogeneous polynomial of degree $d$. The Laurent
series expansion of 
\[ Res_{a_1,\ldots,a_n, a_{n+1} = \ldots = a_{n+r}=
  1} (t^{v_{n+1}} \cdots t^{v_{n+r}} 
P(a_1 t^{v_1},\ldots,a_n t^{v_n}))\]
at the vertex of the secondary
polytope of $\tilde\Delta$ corresponding to the 
triangulation $\tilde\cT$ is
\[ Res_{a_1,\ldots,a_n}(t^{v_{n+1}} \cdots t^{v_{n+r}} P(a_1
t^{v_1},\ldots,a_n t^{v_n})) = 
\sum_{\beta\in R(\overline\Sigma)_\eff} \langle
P(\bar\chi_{1,0},\ldots,\bar\chi_{n,0}) \Phi_\beta
\rangle_{\overline\Sigma_\beta} a^\beta.\] 
\end{theorem} \qed

\section{Mixed residues and mixed volumes}

We keep the notation from the previous section. 

The ring $S_{\tilde\Delta}$ is graded by $\ZZ_{\geq 0}^r$ and
$I_{\tilde\Delta}\subset S_{\tilde\Delta}$ is a homogeneous ideal. For
a partition $k = (k_1,\ldots,k_r)$,
\[ k_1+\ldots+ k_r = n+d, \qquad k_i>0, \]
the restriction of $Res_{a_1,\ldots,a_n}$ to the degree $k$ component
of $I_\Delta$ is called the $k$-mixed residue. 

The following was conjectured by Batyrev and Materov \cite{BM2}:

\begin{theorem} Let $H_f^k$ be the $k$-homogeneous component of
  $H_f$. The $k$-mixed residue of $H_f^k$ is
\[  Res_{a_1,\ldots,a_n} H_f^k = V(\Delta_1^{\bar{k}_1}\cdots
\Delta_r^{\bar{k}_r}),\]
where the right hand side denotes the mixed volume multiplied with
$(n+d-1)!$, and $\bar{k} = (k_1-1,\ldots,k_r-1)$.
\end{theorem}

{\bf Proof.} Let $c_1,\ldots,c_k$ be parameters close to $1$ and
consider the (non-integral) polytope 
\[ \tilde\Delta_c = c_1 \Delta_1 * \ldots * c_r \Delta_r.\]
The volume of $\tilde\Delta_c$ is a polynomial in $c_i$ with
coefficients the normalized mixed volumes:
\[ Vol(\tilde\Delta_c) = \sum_k V(\Delta_1^{\bar{k}_1}\cdots
\Delta_r^{\bar{k}_r}) c_1^{\bar{k}_1}\cdots c_r^{\bar{k}_r}.\]
We may take this as the definition of the mixed volume.

The triangulation $\tilde\cT$ of $\tilde\Delta$ induces a
triangulation $\tilde\cT_c$ of  
$\tilde\Delta_c$ if we replace the vertices $v_i = (\bar{v}_i,e_j)$ by
$v_{i,c} =  (c_j \bar{v}_i,e_j)$, and leave $v_{n+j,c}= v_{n+j} =
(0,e_j)$. If a simplex $\tau\in \tilde\cT$ corresponds to the simplex
  $\tau_c \in \tilde\cT_c$, then an easy determinant computation shows that 
\[ Vol(\tau_c) = Vol(\tau)  c_1^{\bar{k}_1}\cdots c_r^{\bar{k}_r},\]
where $\bar{k}_j = |\{i\in E_j| v_i\in \tau\}|$. 

Let us write $\gamma = (\gamma_1,\ldots,\gamma_{n+r})$, where
$\gamma_i = c_j$ if $i\in E_j$ or if $i=n+j$. We apply
Lemma~\ref{lem-hess} to get:
\begin{align*} H_{f_\gamma} = \sum_k H_f^k c_1^{k_1}\cdots c_r^{k_r}
& \stackrel{RM}{\longmapsto} 
  \sum_{\sigma \in\Sigma} Vol(\sigma) \prod_{v_i\in\sigma} \gamma_i \\
&=  \sum_{\sigma \in \Sigma} Vol(\sigma_c) c_1\cdots c_r \\
&=  Vol(\tilde{\Delta}_c)c_1\cdots c_r \\
&=  \sum_k V(\Delta_1^{\bar{k}_1}\cdots \Delta_r^{\bar{k}_r})
c_1^{k_1}\cdots c_r^{k_r}.
\end{align*}
Comparing the coefficents on both sides we get the desired result. \qed


\begin{thebibliography}{CC}

\bibitem{BBFK} G. Barthel, J.-P. Brasselet, K.-H. Fieseler, L. Kaup,
  {\em Equivariant intersection cohomology of toric varieties},
  Algebraic geometry: Hirzebruch 70 (Warsaw, 1998),  45--68,
  Contemp. Math., 241, Amer. Math. Soc., Providence, RI, 1999.  

\bibitem{BM1} V. V. Batyrev, E. N. Materov, {\em Toric Residues and
  Mirror Symmetry,} preprint {\tt math.AG/0203216}.

\bibitem{BM2} V. V. Batyrev, E. N. Materov, {\em Mixed toric residues and
  Calabi-Yau complete intersections,} preprint {\tt math.AG/0206057}.

\bibitem{B1} L. A. Borisov, {\em Towards the Mirror Symmetry for
  Calabi-Yau Complete intersections in Gorenstein Toric Fano
  Varieties,} preprint {\tt math.AG/9310001}. 

\bibitem{B2} L. A. Borisov, {\em Higher Stanley-Reisner rings and
  toric residues},  preprint {\tt math.AG/0306307}.
 
\bibitem{B} M. Brion, {\em The structure of the polytope algebra},
  Tohoku Math. J. (2)  49  (1997),  no. 1, 1--32. 

\bibitem{BV} M. Brion, M. Vergne, {\em Arrangement of
  hyperplanes. I. Rational functions and Jeffrey-Kirwan residue},
  Ann. Sci. École Norm. Sup. (4)  32  (1999),  no. 5, 715--741. 


\bibitem{CDS} E. Cattani, A. Dickenstein, B. Sturmfels, {\em Residues
  and resultants,}  J. Math. Sci. Univ. Tokyo  5  (1998),  no. 1,
  119--148.

\bibitem{C} D. Cox, {\em Toric residues,}  Ark. Mat.  34  (1996),
  no. 1, 73--96. 

\bibitem{SV} A. Szenes, M. Vergne, {\em Toric reduction and a
  conjecture of Batyrev and Materov}, preprint {\tt math.AT/0306311}.


\end{thebibliography}
\end{document}